\documentclass[letterpaper, 10 pt, conference]{ieeeconf}  

\IEEEoverridecommandlockouts                              

\usepackage{algorithm}
\usepackage[noend]{algpseudocode}
\usepackage{capt-of}
\usepackage{textcomp}
\usepackage{subfigure}
\usepackage{amsmath} 
\usepackage{graphicx}
\usepackage{url}
\usepackage{amssymb,xcolor,trimclip}
\usepackage{hyperref}
\hypersetup{
    colorlinks,
    linkcolor={red!100!black},
    citecolor={blue!50!black},
    urlcolor={blue!80!black}
}

\title{\LARGE \bf
Chance Constraints Integrated MPC Navigation in Uncertainty amongst Dynamic Obstacles: An overlap of Gaussians approach
}

\author{Dhaivat Bhatt$^{1}$ \quad Akash Garg$^{2}$ \quad Bharath Gopalakrishnan$^{1}$   \quad K. Madhava Krishna$^{1}$ 
\thanks{$^{1}$Affiliated with KCIS, Robotics Research Center, IIIT Hyderabad        {\tt\small }}%
\thanks{$^{2}$Affiliated with Delhi Technological University}}%

\begin{document}

\makeatletter
\setlength{\@fptop}{0pt}
\makeatother
\maketitle
\thispagestyle{empty}
\pagestyle{empty}

\newenvironment{rcases}
  {\begin{aligned}}
  {\end{aligned}\rbrace}
  
\small
{
\begin{abstract}
In this paper, we formulate a novel trajectory optimization scheme that takes into consideration the state uncertainty of the robot and obstacle into its collision avoidance routine. The collision avoidance under uncertainty is modeled here as an overlap between two distributions that represent the state of the robot and obstacle respectively. We adopt the minmax procedure to characterize the area of overlap between two Gaussian  distributions, and compare it with the method of Bhattacharyya distance. We provide closed form expressions that can characterize the overlap as a function of control. Our proposed algorithm can avoid overlapping uncertainty distributions in two possible ways. Firstly when a prescribed overlapping area that needs to be avoided is posed as  a confidence contour lower bound, control commands are accordingly realized through a MPC framework such that these bounds are respected. Secondly in tight spaces control commands are computed such that the overlapping distribution respects a prescribed range of overlap characterized by lower and upper bounds of the confidence contours. We test our proposal with extensive set of simulations carried out under various constrained environmental configurations. We show usefulness of proposal under tight spaces where finding control maneuvers with minimal risk behavior becomes an inevitable task. 
\end{abstract}

\section{INTRODUCTION}
Quadcopter MAVs (Micro Aerial Vehicle) are an ideal choice for autonomous
reconnaissance and surveillance because of their small size, high maneuverability, and ability
to fly in very challenging environments. To perform these tasks effectively, the quad-copter
MAV must be able to precisely avoid obstacles while navigating from one point to another. The obstacles include static and dynamic objects as well as other quad-copters operating in the surrounding environment. 
The field of obstacle avoidance for quad-copters has been explored for quite a long time. Many algorithms proposed in past fails to produce desired result because of the uncertainty involved in belief of the MAV/obstacle. A deterministic obstacle avoidance algorithm is not an appropriate in unstructured and uncertain environments. This can lead to substantial degradation in the desired result and can even make the source robot to collide into the obstacle, in the worst case.

In order to deal with the challenges mentioned, we propose a probabilistic Model Predictive Control framework for trajectory optimization of quad-copters. MPC is proven to be an efficient framework due to its receding horizon planning capability. It is an optimization based approach to handle arbitrary number of constraints on state and control.  The framework optimizes the given cost function which takes maneuverability and actuation limitations of the source robot into consideration. In the paper presented, a quadratic goal reaching objective is provided as cost in addition with a jerk cost to obtain a smooth trajectory. We also incorporate actuation constraints to ensure kinematic feasibility of optimal trajectory. This paper has several novel findings and contributes in following ways,

\begin{figure}[t]
\centering
\includegraphics[width=8cm]{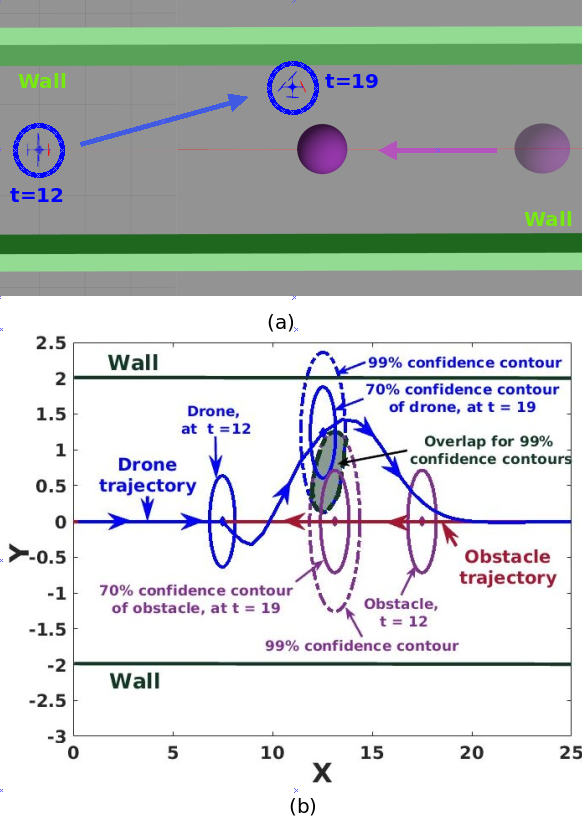}
\caption{We show result of our probabilistic  obstacle avoidance algorithm in constrained corridor when an obstacle is approaching in antipodal configuration. Figure \ref{intro_figure}(a) shows gazebo snapshot of drone positions for two diffrent time instances. In figure \ref{intro_figure}(b). For example, at time t=19, we can see clear overlap between $\mathbf{99\%}$ confidence contours of the drone, however, $\mathbf{70\%}$ confidence contours, which correspond to lower bound are able to avoid penetration.    }
\label{intro_figure}
\end{figure}

\begin{itemize}
   \item This is first such formulation, conditioned on MAV and obstacle uncertainty into an MPC framework through theory of overlapping Gaussians. 
   \item We demonstrate why our modeling is more consistent and appropriate compared to entropic distances between probability distributions
   \item We further introduce this particular uncertainty modeling into an optimal control framework, and jointly optimize in both, the control and overlap parameter space. 
   \item We show effective results in various simulation settings that showcase versatility of the method. Specifically we show where the distributions are non isotropic, which is closer to real setting. 

\end{itemize}

The remainder of this paper is organized as follows. Section \ref{sec_related_work} discusses earlier line of work.  Section \ref{prerequisite} outlines deterministic trajectory optimization framework. Section \ref{COUU} talks about formulating collision avoidance under uncertainty. It also outlines methods to quantify overlap between two Gaussians through various methods. Section \ref{PCAAOBTG} talks about modeling probabilistic collision avoidance using theory built in \ref{TOOG}. In section \ref{TOWCC}, we talk about formulating trajectory optimization algorithm under chance constraints. We show evaluation of trajectory optimization algorithm for a simple case in section \ref{TOR}. In section \ref{Results_section}, we evaluate proposed formulation into an MPC framework for two scenarios. In section \ref{Antipodal_result_section}, we show results for antipodal configuration and in section \ref{OAICC} we show results for tightly bounded constrained corridor setting. In last section, we conclude the paper and discuss future scope of improvement. 


\section{Related work}
\label{sec_related_work}
This section review the recent advances in MPC for autonomous navigation. The evident advantage of using MPC in motion planning and autonomous navigation has been well demonstrated in (\cite{mithun_arxiv},\cite{kim2011control}, \cite{katrakazas2015real}, \cite{schwarting2017parallel}) among many. Formulations along the lines of \cite{mithun_arxiv}, \cite{kim2011control} do a great job in achieving performance in terms of quality of trajectory, computation time and novelty of approach. However, they have been developed for a deterministic setting and hence do not take into consideration the uncertainty in state of the robot and obstacle into their collision avoidance routine. While MPC formulations along the lines of \cite{schwarting2017parallel} does take into consideration the state uncertainty and demonstrate interesting maneuvers in complex driving scenarios, \cite{schwarting2017parallel} considers uncertainty only in the state of the obstacle, and the collision avoidance is modeled through a Minkowski sum approach. Considering robot's uncertainty into Minkowski sum formulation would be very cumbersome as Minkowski sum between two ellipses is very complex. \cite{garimella2017robust} takes into consideration the uncertainty of the drone and assumes the obstacle to be static and determinstic, it models collision as a measure of entropic distance(Similar to Mahalanobis distance\cite{wiki:Mahalanobis_distance}). It is shown in the later section of this paper that formulating probabilistic collision avoidance  as an entropic distance may not be the most appropriate approach when uncertainty in the robot and obstacle is considered. \cite{van2011lqg} attempts to solve for collision avoidance in a multiagent scenario under uncertainty. It achieves it through an RRT framework and a sampling strategy to choose a path corresponding to a desired level of safety confidence. Formulating collision avoidance as a chance constraint is well explored in \cite{gopalakrishnan2017prvo},\cite{gopalakrishnan2015closed}  among the many. \cite{gopalakrishnan2015closed} demonstrates an efficient way of solving an intractable chance constraint through a series of reformulations. These were built on time scale velocity obstacle concepts \cite{singh2013reactive}. In this work, we look at an alternative take on modeling probabilistic collision avoidance as a chance constraint that could seamlessly integrate itself into an MPC framework. We accomplish thus by reformulating collision avoidance as a measure of overlap between two Gaussians(representing state uncertainties of the robot and obstacle). By such a reformulation, we easily avoid the complexity of considering Minkowski sums between two ellipses and other approaches that model probabilistic collision avoidance through entropic distances.

\subsection{Symbols and notations:}
There are several variables of interest to us. We will follow a generic notation scheme. A variable of interest for drone at some time instance $\mathbf{t_i}$ will be expressed in form of $\mathbf{\xi^d_{t_i}}$. For example, uncertainty of drone at time $\mathbf{t_i}$ is denoted by $\mathbf{\Sigma^d_{t_i}}$. A variable of interest for obstacle j at time instance $\mathbf{t_i}$ will be expressed in form of $\mathbf{\xi^{o_j}_{t_i}}$. For example, uncertainty of obstacle j at $\mathbf{t_i}$ is expressed as $\mathbf{\Sigma^{o_j}_{t_i}}$. A variable of mutual interest of obstacle j and drone at time $\mathbf{t_i}$ will be expressed in form of $\mathbf{\xi^{j}_{t_i}}$. For example, overlap between Gaussian populations of obstacle j and drone at time $\mathbf{t_i}$ is defined as $\mathbf{\Upsilon^j_{t_i}}$ in section \ref{TOWCC}.  

\section{Prerequisite} 
\label{prerequisite}

This section describes a deterministic MPC framework, formulated along the lines of \cite{mithun_arxiv}. The motion model considered here is of a holonomic bot, and obstacle avoidance is added as an affine constraint. The entire framework is solved as a sequential convex programming(SCP) routine\cite{scp}. 

\subsection{Trajectory optimization in a deterministic setting}
\label{determinisitcTO}
In deterministic trajectory optimization setting, our objective is to reach the goal in a given amount of time while ensuring a collision-free trajectory. The problem is modeled by considering a set of cost functions and constraints. 

Let the start position of the drone be $\mathbf{X_0}$ = ($\mathbf{x_0}$, $\mathbf{y_0}$, $\mathbf{z_0}$). Our objective is to reach the goal position $\mathbf{G_f}$ =  ($\mathbf{G^x_f}$, $\mathbf{G^y_f}$, $\mathbf{G^z_f}$) in $\mathbf{N}$ time-steps, each time-step of duration $\mathbf{\tau}$. The state of the drone at any time instant $\mathbf{t_i}$ is $\mathbf{X_{t_i}}$ = ($\mathbf{x_{t_i}}$, $\mathbf{y_{t_i}}$, $\mathbf{z_{t_i}}$, $\mathbf{v^x_{t_i}}$, $\mathbf{v^y_{t_i}}$, $\mathbf{v^z_{t_i}}$). Where $\mathbf{V_{t_i}}$ = ($\mathbf{v^x_{t_i}}$, $\mathbf{v^y_{t_i}}$, $\mathbf{v^z_{t_i}}$) is the velocity of the drone at time instant $\mathbf{t_i}$.
We have P obstacles in the environment. Their position at time $\mathbf{t_i}$ is defined as $\mathbf{O^j_{t_i}}$ = ($\mathbf{o^{xj}_{t_i}}$, $\mathbf{o^{yj}_{t_i}}$, $\mathbf{o^{zj}_{t_i}}$), for $\forall \hspace{0.1cm} j = \{1,2,3, ... \hspace{0.1cm} P\}$. For static obstacles, the obstacle locations will be independent of $\mathbf{t_i}$. The drone and obstacles are approximated as circular objects with radius of the drone being $\mathbf{R_{drone}}$ and radius of obstacle $\mathbf{j}$ is $\mathbf{R_{j}}$, $\forall \hspace{0.1cm} j = \{1,2,3,...\hspace{0.1cm} P\}$.

{
\small
\begin{subequations}
\begin{eqnarray}
\label{optimization_model_equations}
\underset{\mathbf{V_{t_i}}}{\mathrm{argmin}} \quad \mathbf{J} = \mathbf{J_{terminal}} +  \mathbf{J_{smooth}} \label{cost_function} \\
\mathbf{X_{t_i + 1}} = f(\mathbf{X_{t_i}}, \mathbf{V_{t_i}}) \label{motion_model} \\
\qquad \mathbf{V_{min}} \leq \mathbf{V_{t_i}} \leq \mathbf{V_{max}} \label{vel_bounds} \\
\qquad \mathbf{a_{min}} \leq \mathbf{\dfrac{V_{{t_i}+1} - V_{{t_i}}}{\tau}} \leq \mathbf{a_{max}} \label{acc_bounds} \\
\qquad \mathbf{C}_{obs_j}(\mathbf{x_{t_i}}, \mathbf{y_{t_i}}, \mathbf{z_{t_i}}, \mathbf{o^{xj}_{t_i}}, \mathbf{o^{yj}_{t_i}}, \mathbf{o^{zj}_{t_i}}, \mathbf{R_{drone}}, \mathbf{R_j}) \leq \mathbf{0} \label{obstacle_avoidance_constraint}
\end{eqnarray}
\end{subequations}
}

The above set of equations defines the cost function as well as constraints. The objective function as described in equation  \ref{cost_function}. 

{
\small
\begin{eqnarray}
\label{terminal_cost}
\mathbf{J_{terminal}} = (\mathbf{x_{t_N} - G_x })^2 + (\mathbf{y_{t_N} - G_y })^2 + (\mathbf{z_{t_N} - G_z })^2
\end{eqnarray}
}

The terminal cost forces our system to achieve goal-state($\mathbf{G_f}$) at the end of the trajectory.

{
\small
\begin{eqnarray}
\label{smoothness_cost}
\mathbf{J_{smooth}} = \sum_{i=2}^{\mathbf{N-1}} (\mathbf{\dfrac{(V_{{t_i}+1} + V_{{t_i}-1} - 2V_{t_i})}{\tau^2}})^2
\end{eqnarray}
}

The smoothness cost as described above ensures smooth trajectory with minimal jerk. It minimizes the jerk which is modeled as second order finite difference between subsequent linear velocities. This term penalizes sudden deviations in the acceleration profile and ensures smooth velocity transitions. 

Equation \ref{motion_model} is the process model of the vehicle. These equations ensure that control variables and states are adhering the motion model of the drone. The motion model for holonomic bot can be described as following,

{
\small
\begin{subequations}
\label{holonomic_motion_model}
\begin{eqnarray}
\mathbf{x_{t_i}} = \mathit{f}(\mathbf{x_0, v^x_{t_1}, v^x_{t_2}, ... v^x_{t_i}, \tau}) = \mathbf{x_0 + \sum_{k=1}^{i} v^x_{t_k}\tau} \\
\mathbf{y_{t_i}} = \mathit{f}(\mathbf{y_0, v^y_{t_1}, v^y_{t_2}, ... v^y_{t_i}, \tau}) = \mathbf{y_0 + \sum_{k=1}^{i} v^y_{t_k}\tau}  \\
\mathbf{z_{t_i}} = \mathit{f}(\mathbf{z_0, v^z_{t_1}, v^z_{t_2}, ... v^z_{t_i}, \tau}) = \mathbf{z_0 + \sum_{k=1}^{i} v^z_{t_k}\tau}
\end{eqnarray}
\end{subequations}
}

Equations \ref{vel_bounds}-\ref{acc_bounds}, represents constraints to model actuation limitations of the drone. The bounds on linear acceleration and velocity ensures that the actuation capabilities of the drone are not violated. 

Equation \ref{obstacle_avoidance_constraint} models collision avoidance constraint between obstacle $\mathbf{j}$ and drone. For a deterministic setting, this will be a simple euclidean distance constraint as below. 

{
\small
\begin{eqnarray}
\begin{aligned}
\label{ob_avoidance}
\qquad \mathbf{C}_{obs_j}(.) = -(\mathbf{x_{t_i}} - \mathbf{o^{xj}_{t_i}})^2-(\mathbf{y_{t_i}} -  \mathbf{o^{yj}_{t_i}})^2 \\  -(\mathbf{z_{t_i}} - \mathbf{o^{zj}_{t_i}})^2 + \mathbf{(R_j + R_{drone})^2} \leq \mathbf{0}.
\end{aligned}
\end{eqnarray}
}

Above constraint is purely non-linear in nature for drone position variable ($\mathbf{x_{t_i}}$, $\mathbf{y_{t_i}}$, $\mathbf{z_{t_i}}$). We linearize it along the lines of \cite{mithun_arxiv} and solve the proposed routine using sequential convex programming. The trajectory optimization routine is then integrated into a model predictive control framework. 

\section{Collision avoidance under uncertainty}
\label{COUU}
Robot motions are generally erroneous in nature. There is always some uncertainty associated with the location of the drone. Apart from that, the sensing module also gives inaccurate estimate of the state of the obstacle. Uncertain position estimate of the obstacle leads to erroneous trajectory estimation. In such cases, the constraint in equation \ref{obstacle_avoidance_constraint} takes the form of \ref{chance_constraint}

{
\small
\begin{eqnarray}
\label{chance_constraint}
\mathbf{\Pr}_j(\mathbf{C}_{obs_j}(\mathbf{x_{t_i}}, \mathbf{y_{t_i}}, \mathbf{z_{t_i}}, \mathbf{o^{xj}_{t_i}}, \mathbf{o^{yj}_{t_i}}, \mathbf{o^{zj}_{t_i}}, \mathbf{R_{drone}}, \mathbf{R_j}) \leq \mathbf{0}) \\
\forall \hspace{0.1cm} j \in \{1,2,3,...\hspace{0.1cm} P\} \nonumber
\end{eqnarray}
}

Constraints of the form \ref{chance_constraint}, are generally known as chance constraints, and in most cases may not have a distribution that can be computed in closed form. The nature of these chance constraints also depends on the form of the deterministic constraints that they are built on. For example chance constraint \ref{chance_constraint} arising out of \ref{ob_avoidance} can be expressed as an entropic distance. The following three sections provide a detailed discussion on three important ways to formulate a chance constraint.

\subsection{Theoretical characterization of chance constraint}
\label{TCOCC}
The chance constraint in \ref{chance_constraint} can take the form of a transformed distribution of \ref{ob_avoidance} as shown in Equation \ref{ob_avoidance},   

{
\small
\begin{eqnarray}
\label{complex_integral}
\idotsint_{V_j} \Pr(\mathbf{D_{t_i}}, \mathbf{O^j_{t_i}}) d\mathbf{D_{t_i}}d\mathbf{O^j_{t_i}}
\end{eqnarray}
}

Where, $\mathbf{D_{t_i}}$ = ($\mathbf{x_{t_i}, y_{t_i}, z_{t_i}}$), position of the drone at time $\mathbf{t_i}$. Under state uncertainty, let $\mathbf{D_{t_i}}$ \texttildelow \hspace{0.1cm} $\mathcal{N}$$\mathbf{(\hat{D}_{t_i}, \Sigma^d_{t_i})}$ and $\mathbf{O^j_{t_i}}$ \texttildelow \hspace{0.1cm} $\mathcal{N}$ $\mathbf{(\hat{O}^j_{t_i}, \Sigma^{o_j}_{t_i})}$ be the Gaussian parameterization of the drone and obstacle $\mathbf{j}$ positions at time $\mathbf{t_i}$. Then, $\mathbf{\Pr(D_{t_i}, O^j_{t_i})}$ takes the following form,

{
\small
\begin{eqnarray}
\label{complex_prb}
\mathbf{\Pr(D_{t_i}, O^j_{t_i})}  \hspace{0.1cm} \thicksim \hspace{0.1cm} \mathcal{N} \mathbf{(}  \begin{pmatrix} \mathbf{\hat{D}_{t_i}}\\ \mathbf{\hat{O}^j_{t_i}}\\  \end{pmatrix},
\left( \begin{array}{cc}
\Sigma^d_{t_i} & 0 \\
0 & \Sigma^{o_j}_{t_i}
\end{array} \right) )
\end{eqnarray}
}

When we substitute equation \ref{complex_prb} in equation \ref{complex_integral}, equation \ref{complex_integral} becomes analytically intractable. The closed form solution of equation \ref{complex_integral} doesn't exist. Authors in \cite{van2011lqg} attempted to tackle problem of multi-robot motion planning for differential drive robots, where the authors numerically evaluate equation \ref{complex_integral} over the region of interest. The region of interest here would be set of positions of the drone and obstacles for which collision occurs. Our objective would be to minimize the value of \ref{complex_integral}, which means we want to maximize the probability collison of avoidance.  However, one  drawback of this procedure is that characterization of such a region($\mathbf{V_j}$) is generally tough. 

Uncertainty matrices($\mathbf{\Sigma^d_{t_i}, \Sigma^{o_j}_{t_i}}$) have been scaled up to accommodate radius values $\mathbf{R_{drone}}$ and $\mathbf{R_j}$. 

There has been a lot of work to characterize  the entropic distance between two distributions. One of the commonly used techniques for entropic distances are chi-square distances, Bhattacharyya distances among the many. We describe case of Bhattacharyya distances which is an extension to Mahalanobis distance\cite{wiki:Mahalanobis_distance}.

\subsection{Bhattacharyya distance}
\label{BCdistance}
Bhattacharyya distance gives measure of similarity between two continuous/discrete probability distributions. It attempts to quantify the overlap between two distributions. For Gaussian distributions, Bhattacharyya distance has well-defined analytical formula to capture the overlap. For two multivariate normal distributions, $\mathbf{p_i}$ = $\mathcal{N}$($\mathbf{\mu_i}$, $\mathbf{\Sigma_i}$) with $\mathbf{i}$ = \{1,2\}, the Bhattacharyya distance metric is defined as following, 

{
\small
\begin{eqnarray}
\label{BC_distance}
\begin{aligned}
\mathcal{BC}(\mathbf{p_1, p_2}) = \frac{1}{8}\mathbf{(\mu_1 - \mu_2)^T\Sigma^{-1}(\mu_1 - \mu_2)} + \\
\frac{1}{2}\ln{\frac{det(\mathbf{\Sigma})}{\sqrt{ det(\mathbf{\Sigma_1})det(\mathbf{\Sigma_2})}}}     
\end{aligned}
\end{eqnarray}
}

Where, $\mathbf{\Sigma}$ = ($\mathbf{\Sigma_1 + \Sigma_2}$)/2. This is an extension of Mahalanobis distance. However, equation \ref{chance_constraint} can not be completely characterized through similarity/dissimilarity given by such entropic measure. Because the notion of probability is not complete in Bhattacharyya metric and particular distance doesn't map to a certain value in probability space. In figure
\ref{BC_explanation},  we demonstrate that for different Gaussian distributions with same overlap have different Bhattacharyya distances. Due to this limitation, entropic distance can't be used to model chance constraint. 

\begin{figure}[h]
\centering
\includegraphics[width=8cm]{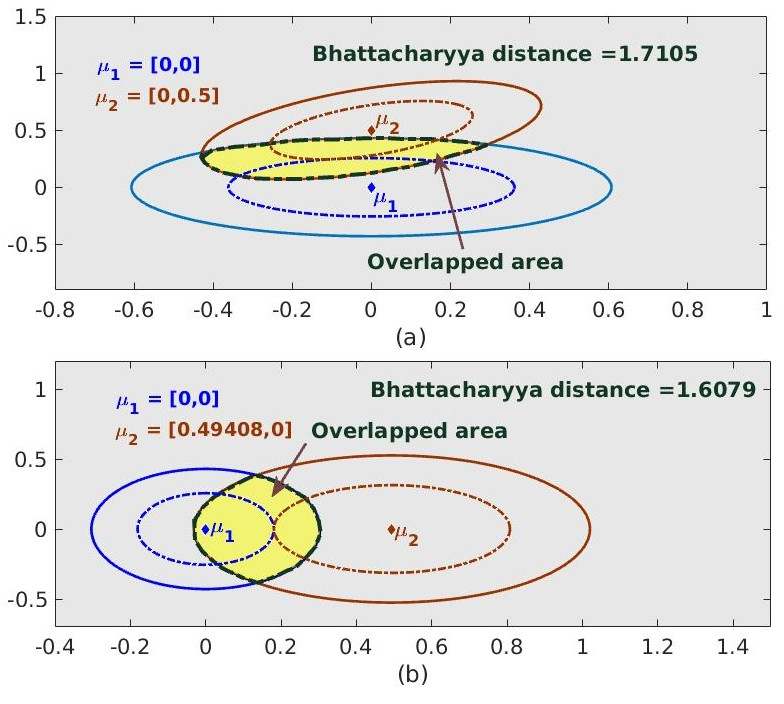}
\caption[Text excluding the matrix]{Diagrammatic explanation of why Bhattacharyya distance is not a perfect metric to model chance constraint problem. In the figure, we have taken $\mathbf{2}$ sets of Gaussian distribution pairs which are touching at same confidence contour of $\mathbf{80.51\%}$. Meaning, area of overlap between two Gaussians shaded in yellow is same for both sets of Gaussian distributions. For figure \ref{BC_explanation}(a), $\mathbf{\Sigma_1} = \left(\protect\begin{smallmatrix}
      0.04 & 0 \\
      0 & 0.02
    \protect\end{smallmatrix}\right)$, 
    $ \mathbf{\Sigma_2} = \left(\protect\begin{smallmatrix}
      0.02 & 0.01 \\
      0.01 & 0.02
 \protect\end{smallmatrix}\right)$. Bhattacharyya distance evaluated using equation \ref{BC_distance} for this set of covariances turns out to be $\mathbf{1.7105}$. While, for figure \ref{BC_explanation}(b),  
 $\mathbf{\Sigma_1} = \left(\protect\begin{smallmatrix}
      0.01 & 0 \\
      0 & 0.02
    \protect\end{smallmatrix}\right)$, 
$ \mathbf{\Sigma_2} = \left(\protect\begin{smallmatrix}
      0.03 & 0 \\
      0 & 0.03    \protect\end{smallmatrix}\right)$. Value of Bhattacharyya distance for figure \ref{BC_explanation}(b) is $\mathbf{1.6079}$. Hence, for same amount of overlap between two sets of Gaussian distributions, Bhattacharyya distances are turning out to be different. Areas shaded with yellow in figure indicate amount of overlap. }
\label{BC_explanation}
\end{figure}


\subsection{Theory of overlapping of Gaussians}
\label{TOOG}
The theory of overlap between two Gaussians has widely been studied in \cite{nowakowska2014tractable}, which is built upon \cite{anderson1962classification}. The authors in \cite{anderson1962classification} attempts to get an optimal linear separator which minimizes the misclassification error when the objective is to classify the sample as coming from one of the several populations. We briefly state the theory to get approximate estimate of component of overlap between two Gaussians. The linear separator proposed in\cite{anderson1962classification} works for two Gaussians of dimension $\mathbf{d  \geq 1}$. 

Let the linear separator(a hyperplane in $\mathbf{d}$ dimensional space) be $\mathbf{\mathbb{\alpha}^Tx = \beta}$ where $\mathbf{\mathbb{\alpha}, x \in \mathbb{R}^d}$ and $\mathbf{\beta \in \mathbb{R}}$. $\mathbf{\mathbb{\alpha}^Tx \leq \beta}$ classifies $\mathbf{x}$ into a first cluster and $\mathbf{\mathbb{\alpha}^Tx > \beta}$ classifies $\mathbf{x}$ into second cluster. We will briefly explain the procedure to obtain $\mathbf{\alpha, \beta}$ and estimate the area of overlap($\mathbf{\Upsilon}$) between two Gaussian distributions $\mathcal{N}$($\mathbf{\mu_i}$, $\mathbf{\Sigma_i}$) with $\mathbf{i}$ = \{1,2\}. 
$\mathbf{x}$ is coming from one of the above two Gaussian distributions. $\mathbf{\alpha^Tx}$ is a transformation which transforms the original distribution into univariate normal distribution. The probability of misclassification when $\mathbf{x}$ is coming from first distribution is,

\setlength\abovedisplayskip{0pt}

{
\small
\begin{align}
\label{P_eta_1}
&\mathbb{P}_1(\mathbf{\alpha^Tx > \beta}) = \mathbb{P}_1(\frac{\mathbf{\alpha^Tx- \alpha^T\mu_1}}{\sqrt{\alpha^T\Sigma_1\alpha}} > \frac{\mathbf{\beta- \alpha^T\mu_1}}{\sqrt{\alpha^T\Sigma_1\alpha}}) \nonumber \nonumber \\
&\mathbb{P}_1(\mathbf{\alpha^Tx > \beta}) = 1 - \mathbf{\Phi}(\frac{\mathbf{\beta - \alpha^T\mu_1}}{\sqrt{\alpha^T\Sigma_1\alpha}}) = 1 - \mathbf{\Phi}(\mathbf{\eta}_1) = \mathbb{P}_1(\mathbf{\eta}_1)
\end{align}
}
\setlength\belowdisplayskip{0pt}

Similarly, probability of misclassification when sample $\mathbf{x}$  belongs to population $\mathbf{2}$ equals, 

\setlength\abovedisplayskip{0pt}
{
\small
\begin{align}
\label{P_eta_2}
&\mathbb{P}_2(\mathbf{\alpha^Tx \leq \beta}) = \mathbb{P}_2(\frac{\mathbf{\alpha^Tx- \alpha^T\mu_2}}{\sqrt{\alpha^T\Sigma_2\alpha}} \leq \frac{\mathbf{\beta- \alpha^T\mu_2}}{\sqrt{\alpha^T\Sigma_2\alpha}}) \nonumber  \\ 
&\mathbb{P}_2(\mathbf{\alpha^Tx \leq \beta}) = 1 - \mathbf{\Phi}(\frac{\mathbf{\alpha^T\mu_2- \beta}}{\sqrt{\alpha^T\Sigma_2\alpha}}) = 1 - \mathbf{\Phi}(\mathbf{\eta}_2) = \mathbb{P}_2(\mathbf{\eta}_2)
\end{align}
}
\setlength\belowdisplayskip{0pt}

Where, $\mathbf{\eta_1 = \frac{\mathbf{\beta - \alpha^T\mu_1}}{\sqrt{\alpha^T\Sigma_1\alpha}}}$ and $\mathbf{\eta_2 = \frac{\mathbf{\alpha^T\mu_2- \beta}}{\sqrt{\alpha^T\Sigma_2\alpha}}}$ are two random variables with univariate standard normal distribution. $\mathbf{\Phi}$ in equations \ref{P_eta_1}-\ref{P_eta_2} denotes a cumulative distribution function for a univariate standard normal distribution. $\mathbf{\Phi}$ is a monotonically increasing function. Our objective is following, 

{
\small
\begin{subequations}
\begin{align}
&max(\mathbb{P}_1(\mathbf{\eta}_1), \mathbb{P}_2(\mathbf{\eta}_2)) \rightarrow 
\underset{\mathbf{\alpha  \in \mathbb{R}^d \beta \in \mathbb{R}}}{\mathrm{min}} \label{max_p_eta} \\
&min(\mathbf{\eta}_1, \mathbf{\eta}_2) \rightarrow 
\underset{\mathbf{\alpha  \in \mathbb{R}^d \beta \in \mathbb{R}}}{\mathrm{max}}
\label{min_eta}
\end{align}
\end{subequations}
}

Equations \ref{max_p_eta}-\ref{min_eta} are equivalent due to monotonic nature of the $\mathbf{\Phi}$. The objective is to maximize the minimum of ($\mathbf{\eta}_1$, $\mathbf{\eta}_2$). The objective is to find $\mathbf{\alpha, \beta}$, which will minimize the maximum probability of misclassification. Analytical characterization of $\mathbf{\alpha, \beta}$ in terms of $\mathbf{\mu_1, \Sigma_1, \mu_2, \Sigma_2}$ can be expressed as following,

{
\small
\begin{subequations}
\begin{align}
\label{Alpha_beta_characterization}
&\mathbf{\alpha} = (\lambda_1\mathbf{\Sigma_1} + \lambda_2\mathbf{\Sigma_2})^{-1}(\mathbf{\mu_2 - \mu_1}) \\
&\mathbf{\beta} = \mathbf{\alpha}^T\mathbf{\mu_1} + \lambda_1 \mathbf{\alpha}^T\mathbf{\Sigma_1}\mathbf{\alpha} = \mathbf{\alpha}^T\mathbf{\mu_2} - \lambda_2 \mathbf{\alpha}^T\mathbf{\Sigma_2}\mathbf{\alpha}
\end{align}
\end{subequations}
}

Here, $\mathbf{\lambda}_1$ and $\mathbf{\lambda}_2$ are two scalars and resulting procedure to estimate these parameters is referred to as \textbf{minmax procedure}. The minmax procedure is an admissible procedure\cite{anderson1962classification} when $\mathbf{\eta}_1$ = $\mathbf{\eta}_2$. For admissible procedure, $\mathbf{\lambda}$ =  $\mathbf{\lambda}_1$ = 1 -  $\mathbf{\lambda}_2$. If we substitute analytical characterization of $\mathbf{\alpha, \beta}$ in $\mathbf{\eta}_1$, $\mathbf{\eta}_2$, the following equality must hold for admissible procedure.

\small{
\begin{align}
\label{base_condition}
\eta_1^2 - \eta_2^2 = \alpha^T[\lambda^2\mathbf{\Sigma_1} - (1 - \lambda)^2\mathbf{\Sigma_2}]\alpha = \mathbf{0}
\end{align}
}

The above criterion is a necessary condition to get the best approximation of the amount of overlap. We will call $\mathbf{\lambda}$ \textit{overlap parameter} as it is a deciding factor which completely characterizes the overlap for a given first and second order Gaussian moments. 


Here, the value of overlap parameter $\mathbf{\lambda}$ is determined heuristically using equation \ref{base_condition} as a base condition. Algorithm to determine $\lambda$ is outlined in \cite{nowakowska2014tractable}.  The overlap parameter $\mathbf{\lambda}$ completely dictates the linear separator parameters $\mathbf{\alpha}$ and $\mathbf{\beta}$.  Once we get optimal value of $\lambda$, we can estimate linear separator parameters $\mathbf{\alpha, \beta}$. Since the procedure is admissible, we can compute $\mathbb{P}_1$($\mathbf{\eta}_1$) = $\mathbb{P}_2$($\mathbf{\eta}_2$) = 
$\mathbb{P}_{minmax}$. The amount of overlap($\mathbf{\Upsilon}$) can be computed as below,

{
\small
\begin{align}
&\mathbf{\eta_1} = \mathit{f}(\alpha, \beta), \quad \alpha = \mathit{f}(\lambda), \quad \beta = \mathit{f}(\lambda) \nonumber \\
\label{overlap_answer}
&\mathbf{\Upsilon} = \mathbb{P}_1(\mathbf{\eta}_1) + \mathbb{P}_2(\mathbf{\eta}_2) = \mathbf{2}\mathbb{P}_{minmax} = \mathit{h}(\mathbf{\mu_1, \mu_2, \Sigma_1, \Sigma_2 , \lambda}) \nonumber  \\
\end{align}
}

Here, we can notice that component of overlap($\mathbf{\Upsilon}$) is parameterized by linear separator parameters $\alpha$ and $\beta$. Overlap $\lambda$ is completely dictating $\alpha$ and $\beta$. Hence, $\lambda$ dictates the amount of overlap $\mathbf{\Upsilon}$. Here, $\lambda$ is determined through an iterative procedure. As seen in figure \ref{Evolution_of_t_explanation}, the overlap is determined by substituting for different values of $\lambda$ in \ref{overlap_answer}, a typical iterative routine settles for some value of $\lambda$, when the overlap is correctly determined, this can be noticed in Fig \ref{Evolution_of_t_explanation}(d)-\ref{Evolution_of_t_explanation}(f).

\begin{figure}[h]
\centering
\includegraphics[width=8cm]{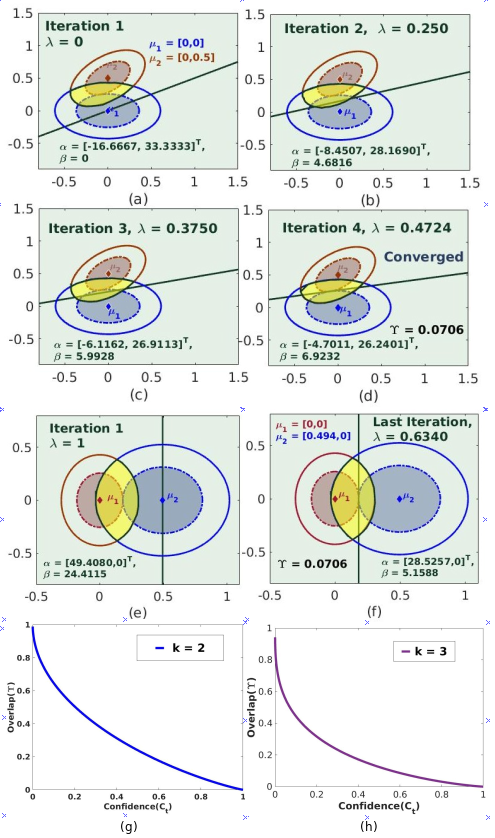}
\caption{This figure demonstrates how the values of overlap parameter $\mathbf{\lambda}$ evolves as algorithm outlined in \cite{nowakowska2014tractable} executes. It starts with initialized value of $\mathbf{\lambda}$ and ultimately converges to the value where the linear separator divides the two Gaussian distributions at the same confidence interval. From figure \ref{Evolution_of_t_explanation}(a) -\ref{Evolution_of_t_explanation}(d), the Gaussian configuration considered is same as that of figure \ref{BC_explanation}(a). Both the Gaussian distributions are touching at confidence contour corresponding to $\mathbf{80.51\%}$. Equality constraint of equation \ref{base_condition} ensures that both the Gaussians are touching at the same confidence interval. The line drawn in the figure is $\mathbf{\alpha^Tx = \beta}$. We can see orientation and position of the line evolving as overlap parameter $\mathbf{\lambda}$ converges. For converged value of $\mathbf{\lambda}$, we can see that the optimal linear separator passes through the point of touch of the two Gaussian distributions. Figure \ref{Evolution_of_t_explanation}(e) -\ref{Evolution_of_t_explanation}(f) shows $\mathbf{\lambda}$ converging for Gaussian configuration considered in figure \ref{BC_explanation}(b). We can see value of overlap turning out to be same for both sets as shown in figure \ref{Evolution_of_t_explanation}(d)-\ref{Evolution_of_t_explanation}(f).  } 
\label{Evolution_of_t_explanation}
\end{figure}


\begin{figure*}[t]
\centering     
\subfigure[]
{\label{fig:a}\includegraphics[width=57mm]{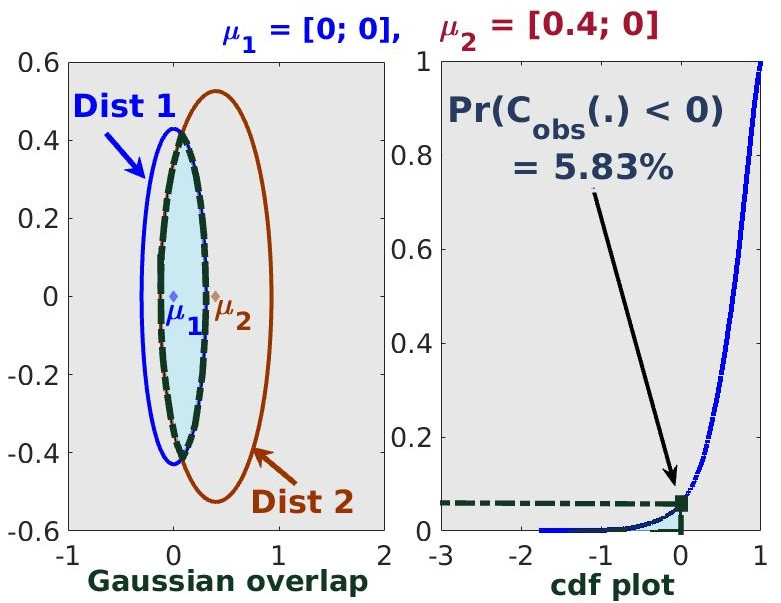}}
\subfigure[]{\label{fig:b}\includegraphics[width=57mm]{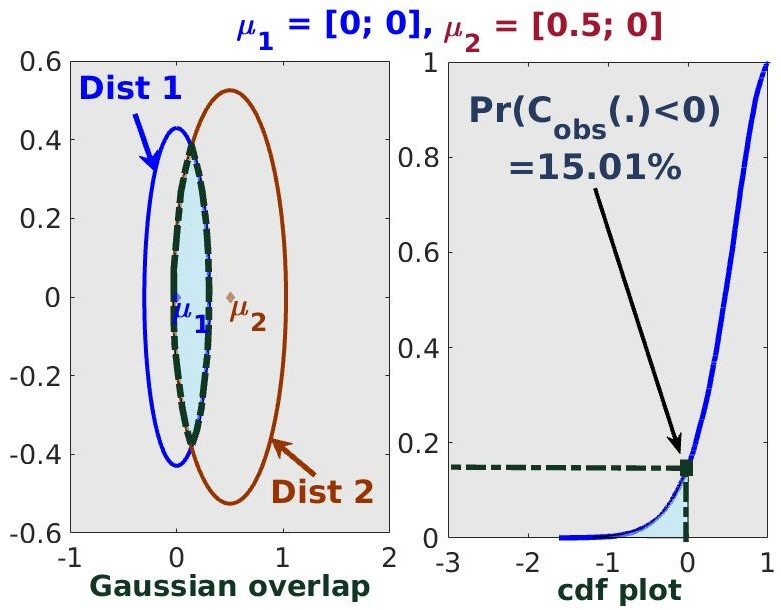}}
\subfigure[]
{\label{fig:c}\includegraphics[width=57mm]{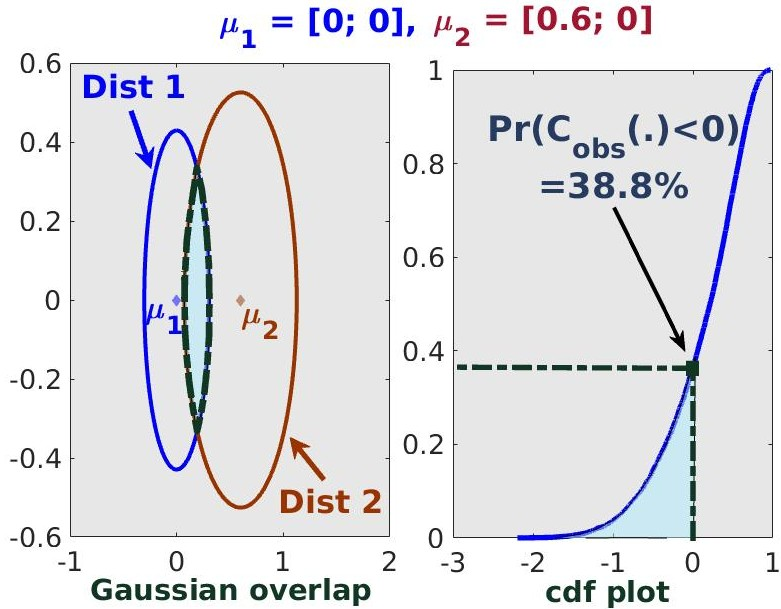}}
\caption[Text excluding the matrix]{We establish an analogy between minimization of overlap between Gaussians as maximization of probability of collision avoidance modeled according to equation \ref{chance_constraint}. Through an illustration, we demonstrate that as the amount of overlap between two distributions decreases, probability of collision avoidance increases. This is aptly conveyed through figure \ref{fig:a}-\ref{fig:b}-\ref{fig:c}. For example, in figure \ref{fig:a}, for a significant amount of overlap, the area below $\mathbf{0}$ in the cdf plot is very less. However, in figure \ref{fig:b} and \ref{fig:c}, as the overlap decreases(achieved through incrementally separating $\mathbf{\mu_2}$ from $\mathbf{\mu_1}$), the corresponding area below $\mathbf{0}$ observed in the cdf plots significantly increases. Thus, it is clear that as the overlap between these two distributions decreases, the probability of collision avoidance[\ref{chance_constraint}](conveyed through cdfs) increases. The cdf plots of equation \ref{chance_constraint} were generated through \textit{ecdf()} function of Matlab.
 The covariances considered for this demonstration are, $\mathbf{\Sigma_1} = \left(\protect\begin{smallmatrix}
      0.01 & 0 \\
      0 & 0.02
   \protect\end{smallmatrix}\right)$ and 
$ \mathbf{\Sigma_2} = \left(\protect\begin{smallmatrix}
      0.03 & 0 \\
      0 & 0.03    \protect\end{smallmatrix}\right)$.  }  \label{chance_constraint_intuition}
\end{figure*}

When we estimate area of overlap for the two sets of Gaussian distributions in figure \ref{BC_explanation} using algorithm outlined in \cite{nowakowska2014tractable}, the area of overlap turns out to be $\mathbf{\Upsilon}$ = 0.0706, which is same for both sets of Gaussian distributions shown in figures \ref{BC_explanation}(a) and \ref{BC_explanation}(b). This leads to the fact that for any two 2D Gaussian distributions touching at the confidence interval of $\mathbf{80.51\%}$, their area of overlap will always be equal to 0.0706. In other words, there is a unique mapping between area of overlap of two $\mathbf{k}$ variate Gaussians and confidence contours $\mathbf{c_t}$(where they are touching each other). The current scope of this paper explores this concept for $\mathbf{k = 2,3}$. We create a table that would give us unique value of overlap for a given value of $\mathbf{c_t}$. We show relation between contour of touch($\mathbf{c_t}$) and area of overlap($\mathbf{\Upsilon}$), in the figure \ref{Evolution_of_t_explanation}(g)-\ref{Evolution_of_t_explanation}(h).

\section{Probabilistic collision avoidance as overlap between two Gaussians}
\label{PCAAOBTG}
Equation \ref{complex_integral} represents analytical expression of chance constraint defined in section \ref{COUU}. Our objective is to maximize the probability of collision avoidance. Here, we propose a novel formulation by posing chance constraint problem as minimization of overlap between component of two Gaussian distributions. We demonstrate in figure \ref{chance_constraint_intuition} that minimization of overlap between two Gaussians is analogous to maximization of probability of collision avoidance.
$\newline$
\subsection{Trajectory optimization with chance constraints}
\label{TOWCC}
In section \ref{prerequisite}, we discussed deterministic  trajectory optimization routine. In this section, we reformulate this routine to accommodate state uncertainty. We express collision avoidance constraint as desired measure of overlap between Gaussian populations of drone and obstacle. The remodeling will use theory of overlapping of Gaussians described in section \ref{TOOG}. Since there is an uncertainty in drone's position, we will redefine $\mathbf{J_{terminal}}$ as below,

\small
{
\begin{align}
\label{JTermNew}
\mathbf{J_{terminal}} = \mathbf{(D_{t_N} - G_f)(\Sigma^d_{t_N})^{-1} (D_{t_N} - G_f)^T}
\end{align}
}

Equation \ref{JTermNew} is mahalanobis distance\cite{wiki:Mahalanobis_distance}, which characterizes the number of standard deviations a point is away from mean of a distribution. Equation \ref{JTermNew} will minimize number of standard deviations goal point is away from the mean position of the drone at the end of the trajectory. 

Here, we are assuming no uncertainty in the actuation. Further, we assume that belief propagation for drone and obstacles for a given time-horizon is known. Hence, minimization of overlap between drone and obstacle populations can be thought of as minimum number of standard deviations$(\mathbf{c_{min}})$ a drone should deviate from its path to ensure collision free trajectory. Let that number be denoted by $\mathbf{c_{t_i}}$ at time instance $\mathbf{t_i}$. While planning the trajectory, our constraint is to ensure that the minimum value of 
$\mathbf{c_{t_i}}$ for $\mathbf{i \in \{1,2..., N\}}$ is larger than certain threshold $\mathbf{c_{min}}$. Which is analogous to saying that the overlap between drone and obstacle at any time instant should not be greater than overlap threshold $\mathbf{\Upsilon_{max}}$. 

\subsubsection{Reformulation of collision avoidance constraint}

Our goal is to reach $\mathbf{G_f}$ from the start position $\mathbf{X_0}$ in $\mathbf{N}$ time-steps, each time-step of duration $\mathbf{\tau}$. Here, our objective is to find optimal set of velocity commands $\mathbf{V_{t_i}}$ = ($\mathbf{v^x_{t_i}}$, $\mathbf{v^y_{t_i}}$, $\mathbf{v^z_{t_i}}$) which would satisfy our constraints as well as minimize the cost. We will be using the process model of the drone explained in section \ref{prerequisite}.

Let the overlap between drone and obstacle $\mathbf{j}$ at time instance $\mathbf{t_i}$ be $\mathbf{\Upsilon^j_{t_i}}$, which is  dictated by overlap parameter $\mathbf{\lambda^j_{t_i}}$. then,

{
\small
\begin{subequations}
\begin{align}
&\mathbf{D_{t_i}} = \mathit{f}(\mathbf{V_{t_1}, V_{t_2},V_{t_3} .... V_{t_i}}) \\
&\mathbf{\Upsilon^j_{t_i}} = \mathit{g_1}(\mathbf{D_{t_i}, \Sigma^d_{t_i}, \hat{O}^j_{t_i}, \Sigma^{o_j}_{t_i}, \lambda^j_{t_i}}) \label{overlap_characterization} \end{align}
\end{subequations}
}

In equation \ref{overlap_characterization}, overlap is expressed in terms of drone/obstacle positions and their corresponding uncertainties. Drone position is expressed in terms of control commands($\mathbf{V_{t_1}, V_{t_2}, ... V_{t_i}}$) as explained in motion model(equation \ref{holonomic_motion_model}).  Hence, equation \ref{overlap_characterization} is completely parameterized by control commands($\mathbf{V_{t_1}, V_{t_2}, ... V_{t_i}}$) and overlap parameter($\mathbf{\lambda^j_{t_i}}$). We can express condition to admissible procedure(equation \ref{condition_to_admissible_procedure}) in terms of control and overlap parameter, 

{
\small
\begin{subequations}
\begin{align} 
\begin{split}
&(\mathbf{\eta^{d}_{t_i}})^2 - (\mathbf{\eta^{o_j}_{t_i}})^2 = \mathit{g_2}(\mathbf{D_{t_i}, \Sigma^d_{t_i}, \hat{O}^j_{t_i}, \Sigma^{o_j}_{t_i}, \lambda^j_{t_i}})  \\  \label{condition_to_admissible_procedure}
&(\mathbf{\eta^{d}_{t_i}})^2 - (\mathbf{\eta^{o_j}_{t_i}})^2= \mathit{f_2}(\lambda^j_{t_i}, \mathbf{V_{t_1}}, \mathbf{V_{t_2}}, \mathbf{V_{t_3}} ...  \mathbf{V_{t_i}})
\end{split} \quad \bigg \rbrace \text{} \\
&\mathbf{\Upsilon^j_{t_i}} = 2\mathbb{P}^d_{t_i}\mathbf{(\eta^d_{t_i})} = 2\mathbb{P}^{o_j}_{t_i}\mathbf{(\eta^{o_j}_{t_i})}  = \mathit{f_1}(\lambda^j_{t_i},\mathbf{V_{t_1}}, \mathbf{V_{t_2}}, \mathbf{V_{t_3}} ...  \mathbf{V_{t_i}}) \label{overlap_final_value}
\end{align}
\end{subequations}
}

Equation \ref{condition_to_admissible_procedure}-\ref{overlap_final_value} are in accordance with equation \ref{base_condition}-\ref{overlap_answer}, equation \ref{condition_to_admissible_procedure} models necessary condition for the procedure to be admissible in nature. 

The maximum allowed overlap $\mathbf{\Upsilon_{max}}$ is uniquely related to minimum number of standard deviations($\mathbf{c_{min}}$) a drone should deviate in order to avoid obstacle with certain minimum confidence. $\mathbf{c_{min}}$ is expressed in terms of confidence intervals directly. So for a particular confidence interval, $\mathbf{c_{min}}$ will have a unique scalar value. Hence, a unique $\mathbf{\Upsilon_{max}}$ value as explained in section \ref{TOOG}. So, a chance constraint in terms of overlap between two Gaussians can be expressed as below. 
{
\small
\begin{eqnarray}
{C}_{obs_j}(.)  = 
\begin{cases}
    \mathbf{\Upsilon^j_{t_i} =  \mathit{f^{lin}_1}(\lambda^j_{t_i}, \mathbf{V_{t_i}})  \leq \Upsilon_{max}}  \\
    (\mathbf{\eta^{d}_{t_i}})^2 - (\mathbf{\eta^{o_j}_{t_i}})^2 = \mathit{f^{lin}_2}(\lambda^j_{t_i}, \mathbf{V_{t_i}}) = \mathbf{0} 
\end{cases}
\label{new_obs_avoidance_constraint}
\end{eqnarray}
}
Both sub-constraints of collision avoidance constraint(equation \ref{new_obs_avoidance_constraint}) are parameterized by velocity controls($\mathbf{V_{t_1}, V_{t_2}, ... V_{t_i}}$) and overlap parameter($\mathbf{\lambda^j_{t_i}}$). $\mathbf{\Upsilon^j_{t_i} \leq \Upsilon_{max}}$ ensures that the drone is avoiding the obstacle with certain minimum confidence. $(\mathbf{\eta^{d}_{t_i}})^2 - (\mathbf{\eta^{o_j}_{t_i}})^2 = 0$ ensures that the overlap parameter($\mathbf{\lambda^j_{t_i}}$) we get through SCP routine is optimal and satisfies condition to admissible procedure. Closed form expressions for $\mathit{f_1}$(.) and $\mathit{f_2}$(.) are functions of optimization variables($\mathbf{V_{t_i}}$, $\mathbf{\lambda^j_{t_i}}$) which are computed using mathematica\cite{mathematica}. They are non-linear in our variables of interest. We linearize them as shown below,

{
\small
\begin{align}
&\mathit{f^{lin}_1()} = \mathit{\bar{f}_1()} +  \sum_{k=1}^{i}\triangledown_{\mathbf{V_{t_k}}}(\mathbf{V_{t_k} - \bar{V}_{t_k}}) +  \triangledown_{\mathbf{\lambda^j_{t_i}}}(\mathbf{\lambda^j_{t_i} - \bar{\lambda}^j_{t_i}}) \label{f1lin} \\
&\mathit{f^{lin}_2()} = \mathit{\bar{f}_2()} +  \sum_{k=1}^{i}\triangledown_{\mathbf{V_{t_k}}}(\mathbf{V_{t_k} - \bar{V}_{t_k}}) +  \triangledown_{\mathbf{\lambda^j_{t_i}}}(\mathbf{\lambda^j_{t_i} - \bar{\lambda}^j_{t_i}}) \label{f1lin} 
\end{align}
}

$\mathit{f^{lin}_1}$(.) and $\mathit{f^{lin}_2}$(.) are affine approximations of $\mathit{f_1(.)}$ and $\mathit{f_2(.)}$. $\triangledown_{\mathbf{V_{t_k}}}$ and $\triangledown_{\mathbf{\lambda^j_{t_i}}}$ are partial derivatives with respect to $\mathbf{V_{t_k}}$ and $\mathbf{\lambda^j_{t_i}}$ respectively. 

\subsubsection{Trajectory optimization algorithm}
We outline complete trajectory optimization algorithm built on above scheme. Let the trajectory of obstacle $\mathbf{j}$ be denoted by $\mathbf{\mathit{\Omega_j}}$ = \{$\mathbf{O^j_{t_1}, O^j_{t_2}, O^j_{t_3}, ... O^j_{t_N}}$\} and trajectory of all $\mathbf{P}$ obstacles be $\mathbf{\Pi_P}$ = \{$\mathit{\Omega_1, \Omega_2, \Omega_3, ... \Omega_P}$\}. Let overlap parameters between obstacle $\mathbf{j}$ and drone for $\mathbf{N}$ timesteps be $\mathbf{\Lambda_j}$ = \{$\mathbf{\lambda^j_{t_1}, \lambda^j_{t_2} ... \lambda^j_{t_N}}$\}. Algorithm
$\mathbf{1}$ outlines proposed SCP routine where we jointly optimize over control($\mathbf{V_{t_i}}$) and overlap parameter($\mathbf{\lambda^j_{t_i}}$) space.

\begin{algorithm}
\caption{ProbabilisticTrajOpt($\mathbf{\Upsilon_{max}}$, $\mathbf{\Pi_P}$, $\mathbf{\Sigma_{drone}}$, $\mathbf{\Sigma_{obstacle}}$)}
\begin{algorithmic}[1]
\State \textbf{Initialization:} Guess for $\mathbf{\bar{\Lambda}^k_{j}(t)}$, $\mathbf{\bar{V}^k(t)}$, iteration counter $\mathit{k}$ = 0 
\State $\mathbf{\bar{D}^\mathit{k}(t)}$ = $\mathit{InitializeTrajectory}$($\mathbf{\bar{V}^\mathit{k}(t)}$)
\State \textbf{while} $\mid$ $\mathbf{J_{\mathit{k}+1} - J_{\mathit{k}}}$ $\mid$ $\mathbf{\geq \delta}$ \textbf{do}
\begin{align*}
&\mathbf{V^\mathit{k}(t), \Lambda^\mathit{k}_{j}(t)} = \mathrm{argmin} \quad \mathbf{J_\mathit{k}} \\
&\textbf{subject to} \\
&\qquad \mathbf{X_{t_i + 1}} = \mathit{f}(\mathbf{X_{t_i}, V_{t_i}}) \\
&\qquad  \mathbf{V_{min}} \leq \mathbf{V_{t_i}} \leq \mathbf{V_{max}} \\
&\qquad  \mathbf{a_{min}} \leq (\mathbf{\dfrac{V_{{t_i}+1} - V_{{t_i}}}{\tau}}) \leq \mathbf{a_{max}} \\
& \qquad {C}_{obs_j}(\mathbf{\bar{ \Lambda}^\mathit{k}_{j}(t)}, \mathbf{\bar{V}^\mathit{k}(t)}) \mathbf{ \leq 0}, \hspace{0.1cm} \forall \hspace{0.1cm} j = \{1,2,3,...\hspace{0.1cm} P\} 
\end{align*}
\State \quad $\mathbf{\bar{\Lambda}^\mathit{k}_{j}(t)} \leftarrow \mathbf{\Lambda^\mathit{k}_{j}(t)}$ 
\State \quad $\mathbf{\bar{V}^\mathit{k}(t)} \leftarrow \mathbf{V^\mathit{k}(t)}$ 
\State \quad $\mathit{k} \leftarrow \mathit{k+1}$
\State \textbf{end while}
\end{algorithmic}
\end{algorithm}
}

\subsection{Evaluation of trajectory optimization}
\label{TOR}
In this section, we evaluate proposed trajectory optimization routine for single obstacle-drone configuration. The start position of the drone is $\mathbf{X_0 = [0,0,0]}$ and it has to reach destination $\mathbf{G_f = [10,0,0]}$ in $\mathbf{N=20}$ timesteps of duration $\mathbf{\tau = 1.0}$ seconds. Obstacle starts from $\mathbf{O^1_{t_1} = [10,0,0]}$ to reach $\mathbf{O^1_{t_N} = [0,0,0]}$. Our objective is to find an optimal trajectory where $\mathbf{c_{min} = 60\%}$ confidence contour of the drone avoids $\mathbf{60\%}$ confidence contour of the obstacle. Meaning, during the entire trajectory, $\mathbf{60\%}$ confidence contour of drone should not penetrate $\mathbf{60\%}$ confidence contour of obstacle. When two $\mathbf{3}$ dimension Gaussians touch each other at $\mathbf{60\%}$ confidence contours, the area of overlap is $\mathbf{0.0861}$. Overlap between drone and obstacle populations at any point of time during the trajectory can not be more than $\mathbf{\Upsilon_{max} = 0.0861}$ . We use algorithm $\mathbf{1}$ to get optimal trajectory satisfying the constraints.

\begin{figure}[h]
\centering
\includegraphics[width=8cm]{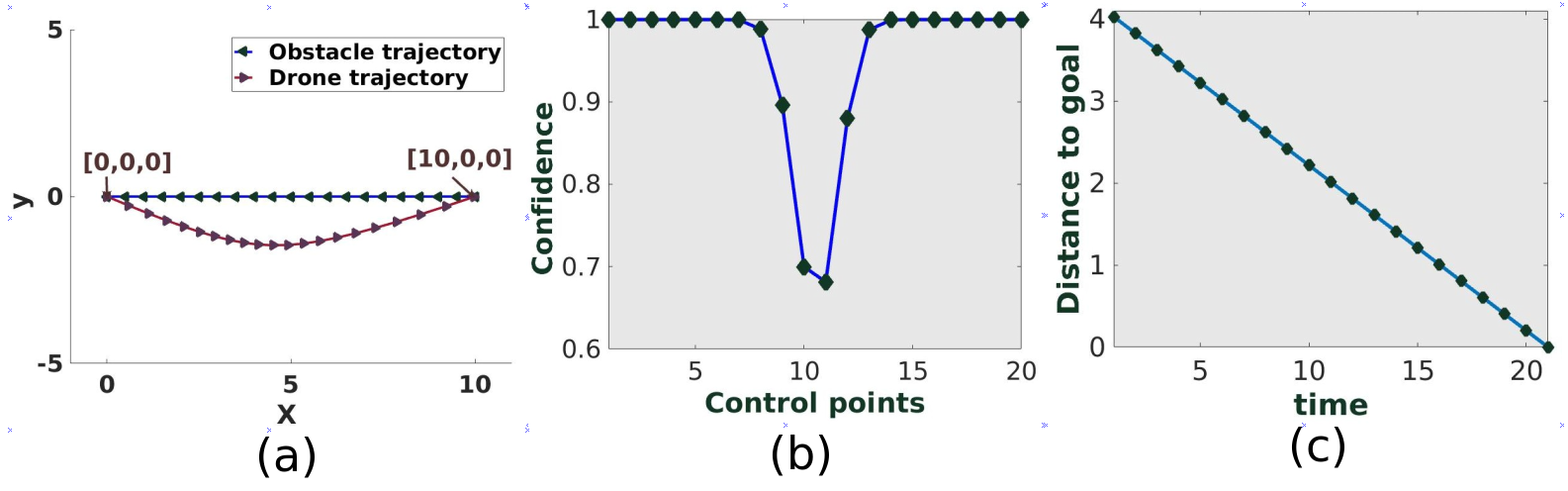}
\caption{In this figure, we show evaluation of trajectory optimization routine proposed in Algorithm $\mathbf{1}$. In figure \ref{Traj_opt_demo}(a), we see trajectory of drone and obstacle for a simple scenario. Figure \ref{Traj_opt_demo}(b) shows confidence contour avoidance profile for this case. We see that as drone and Obstacle keep coming closer, their confidence contours corresponding to confidence = $\mathbf{60\%}$ avoid each-other. \ref{Traj_opt_demo}(c) shows how Mahalanobis distance\cite{wiki:Mahalanobis_distance} between goal point and drone reduces over time. }
\label{Traj_opt_demo}
\end{figure}

\section{Results and discussions}
\label{Results_section}
We construct a model predictive control framework by using proposed probabilistic trajectory optimization routine as a base. Our proposal has been extensively evaluated for wide range of safety critical configurations. We show two such challenging situations to evaluate our proposal. We use Matlab based CVX\cite{cvx} to prototype many of these scenarios. For a faster implementation, we use python based CVXOPT\cite{cvxopt}. The simulations are carried out using Rotors\cite{furrer2016rotors}, which is a micro aerial vehicle simulation framework built in gazebo\cite{koenig2004design}. The proposed approach is implemented in a model predictive control framework(MPC) along the lines of \cite{raemaekers2007design}, \cite{raemaekers2007design} uses the idea of a receding horizon as the basis for building an MPC. We plan for a finite horizon and upto some intermediate way-point along the original trajectory. Intermediate way-points are placed at regular intervals to avoid extreme deviations in drone trajectory. In this section, we show results for two interesting applications. The detailed video of simulations under various safety critical situations is provided at \url{http://robotics.iiit.ac.in/people/dhaivat.bhatt/CDC_video/index.html}

\subsection{Antipodal configuration}
\label{Antipodal_result_section}
\begin{figure*}[h]
\centering     
\includegraphics[width=16cm]{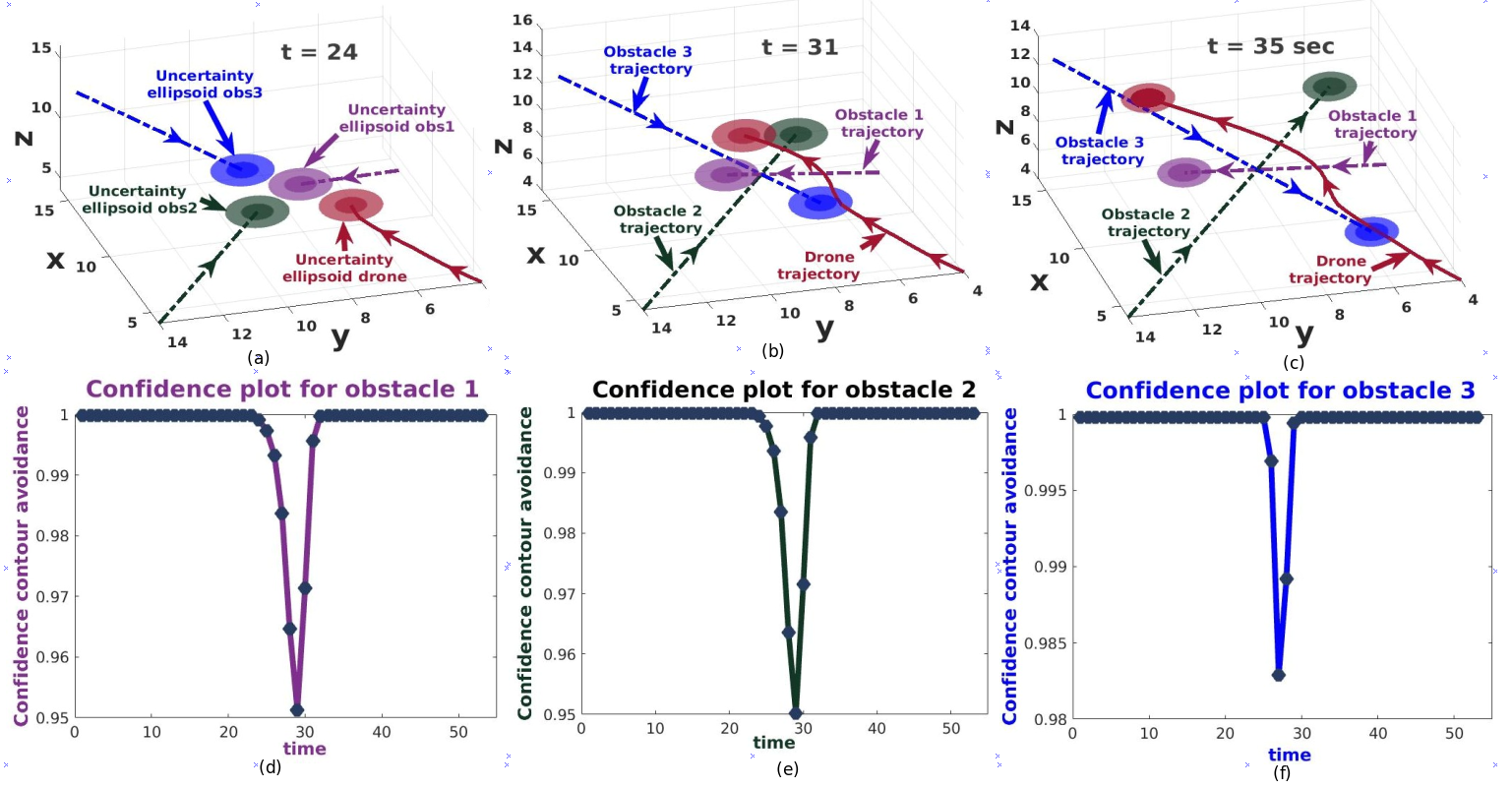}
\caption{Antipodal setting: The drone adopts a maneuver relating to a particular confidence of safety when it encounters obstacles in its sensing range. This is clearly shown in figure \ref{antipodal_result}(a)-\ref{antipodal_result}(c). Figure \ref{antipodal_result}(a), shows the situation, where the drone encounters obstacles, withing its sensor range and starts taking appropriate control actions.Figure \ref{antipodal_result}(b), highlights the resultant maneuver, that the drone adopts to achieve a targeted level of safety. Figure \ref{antipodal_result}(c) shows the goal reaching ability of the drone, after avoiding obstacles. The lighter ellipsoidal shades in these figures represent the uncertainty region encompassing the mean positions of the drone and obstacles(filled with darker shades). Figures \ref{antipodal_result}(d)-\ref{antipodal_result}(e)-\ref{antipodal_result}(f), shows the plots of confidence intervals for the collision avoidance maneuvers the that drone adopted in figures \ref{antipodal_result}(a)-\ref{antipodal_result}(c). Our constraint was to ensure that $\mathbf{90\%}$ confidence contours of drone avoids at least $\mathbf{90\%}$ confidence contours of all obstacles. From figures \ref{antipodal_result}(d)-\ref{antipodal_result}(e)-\ref{antipodal_result}(f), we can observe that maximum overlap between the drone and any obstacle is $\mathbf{5\%}$,, i.e. only the $\mathbf{95\%}$ confidence contours grace each-other. }
\label{antipodal_result}
\end{figure*}

In this experiment, we show 3 obstacles attacking drone in an antipodal configuration. The drone detects obstacles at sensing range of $\mathbf{S_r = 10}$ meters. The acceleration bounds $\mathbf{a_{min}, a_{max}}$ are $\mathbf{-0.5}$ $m/s^2$, $\mathbf{0.5}$ $m/s^2$ respectively. Value of $\mathbf{V_{min}}$ and $\mathbf{V_{max}}$ are $\mathbf{0}$ $m/s$ , $\mathbf{3.0}$ $m/s$ respectively. We keep a planning horizon of $\mathbf{28}$ steps, each of duration $\mathbf{\tau = 0.3}$ seconds. We keep re-planning at every $\mathbf{0.3}$ seconds. Radius of drone/obstacles is taken as $\mathbf{0.5}$ meters. For this configuration, we consider $\mathbf{c_{min} = 90\%}$. In other words, our objective is to ensure that at least $\mathbf{90\%}$ confidence contours of drone and obstacles do not penetrate each other during entire journey of the drone. As soon as it detects the obstacles, it starts deviating from its trajectory and $\mathbf{90\%}$ confidence contour of the drone avoids $\mathbf{90\%}$ confidence contours of all 3 obstacles. The position uncertainties considered are  $ \mathbf{\Sigma_{drone}} = \mathbf{\Sigma_{obs1}} = \mathbf{\Sigma_{obs2}} =
\mathbf{\Sigma_{obs3}} =
\left(\protect\begin{smallmatrix}
      0.02 & 0 & 0 \\
      0 & 0.02 & 0 \\
      0 & 0 & 0.02 \\
      \protect\end{smallmatrix}\right)$. With this configuration, we show comprehensive results in figure \ref{antipodal_result}. Sequence of images in figure \ref{antipodal_result} shows some snaps during obstacle avoidance maneuver, the navigation was successful and the lower bound was respected throughout the trajectory. As shown in confidence plots in figure \ref{antipodal_result}(d)-\ref{antipodal_result}(e)-\ref{antipodal_result}(f),  the maximum penetration for violet and dark green obstacles was at $\mathbf{5\%}$, while for blue obstacle, it was $\mathbf{2\%}$.  

\vspace{0.7cm}

\begin{figure*}[h]
\centering     
\includegraphics[width=16cm]{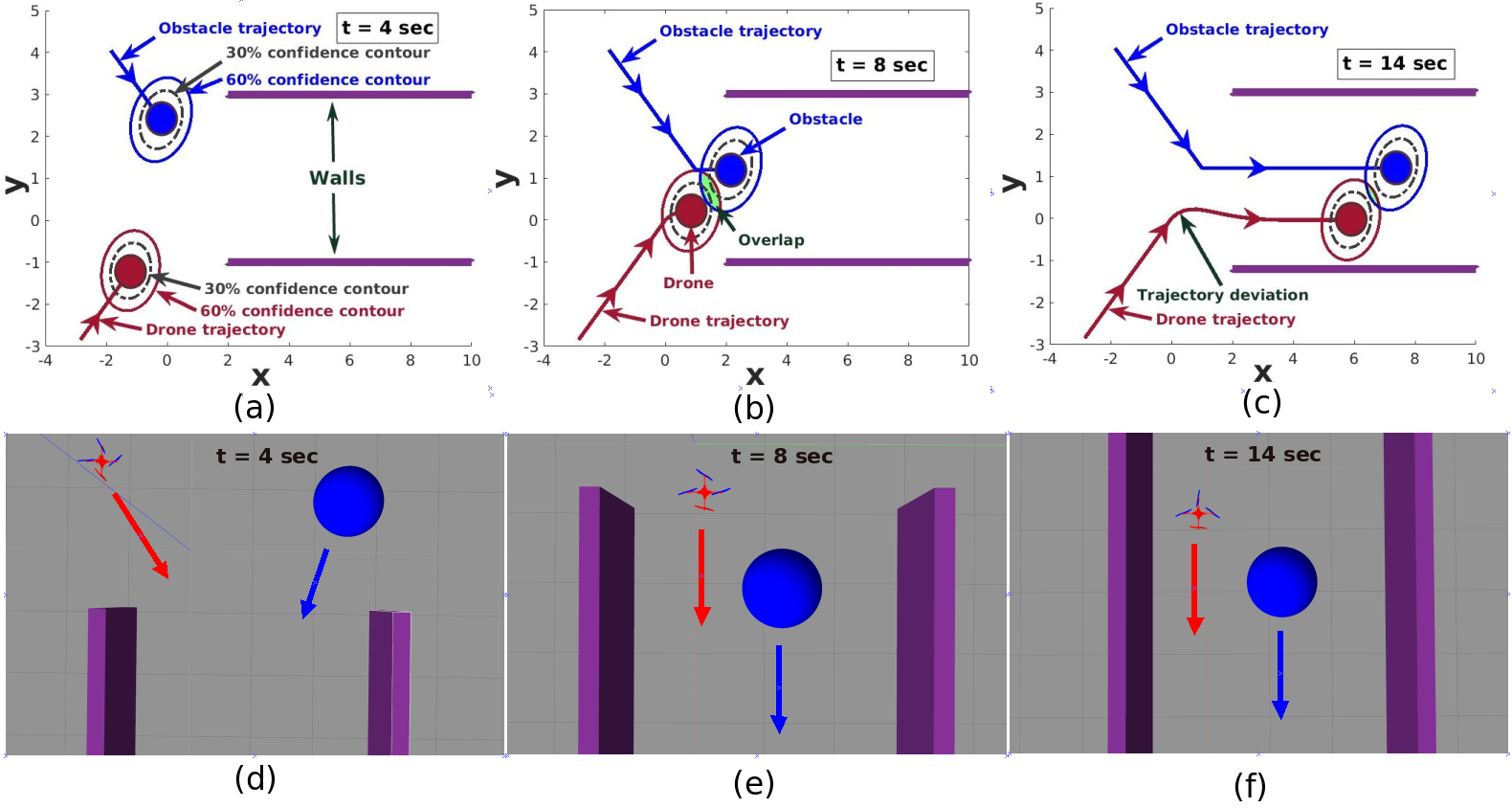}
\caption{Constrained corridor setting: In figure \ref{T_traj_results}(a), confidence contour for 30\% and 60\% are shown for both, drone and obstacle. Drone senses presence of obstacle at time $\mathbf{t = 4 sec}$. As drone enters in corridor, lower and upper bound constraints are enforced. Which ensures minimum risk behavior. For example, in figure \ref{T_traj_results}(b)- \ref{T_traj_results}(c), we can see overlap between drone and obstacle shaded in green for 60\% confidence contours. While, 30\% confidence contours which correspond to lower bound never penetrate each other as visible in figure \ref{T_traj_results}(b)-\ref{T_traj_results}(c). Figure \ref{T_traj_results}(e) shows top view of figure \ref{T_traj_results}(b) in gazebo. Under such tightly bounded spaces, we observe that drone is able to safely maneuver constraints without crashing into walls. Trajectory deviation shown in \ref{T_traj_results}(c) depicts that drone can maneuver in a way that would stay reasonably behind the obstacles while satisfying both lower and upper bound constraints. Figures \ref{T_traj_results}(d)-\ref{T_traj_results}(e)-\ref{T_traj_results}(f) are gazebo results of figures \ref{T_traj_results}(a)-\ref{T_traj_results}(b)-\ref{T_traj_results}(c) in bird's eye view. }
\label{T_traj_results}
\end{figure*}

\subsection{Obstacle avoidance in constrained corridor} 
\label{OAICC}
\begin{figure}[h]
\centering     
\includegraphics[width=8cm]{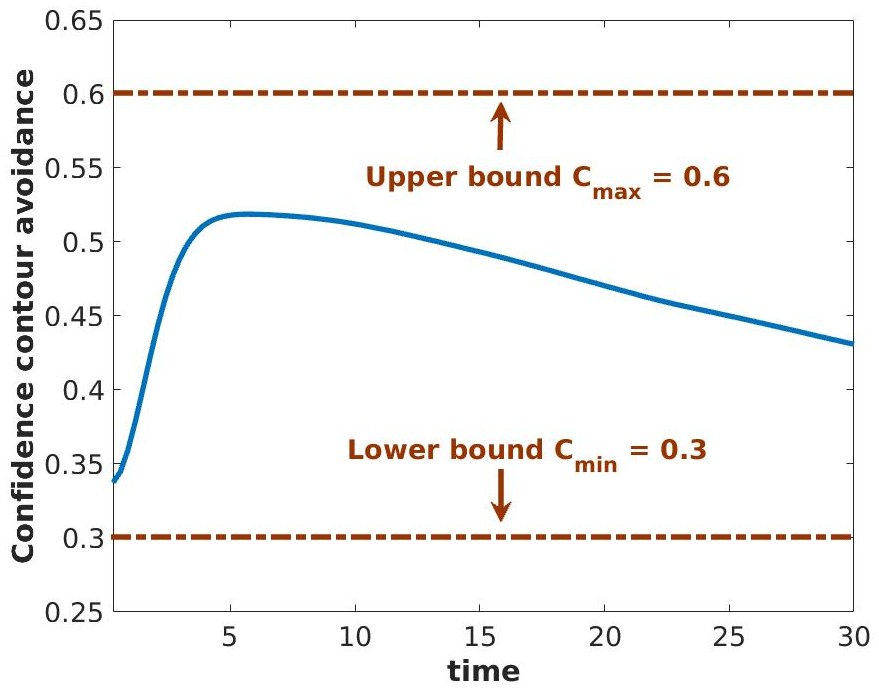}
\caption{This figure shows how confidence contour of touch($\mathbf{c_{t_i}}$) is changing over time for setting considered in \ref{OAICC}.  Throughout the trajectory, our drone is able to satisfy  both, lower and upper-bound constraints. Initial rise in plot suggests that drone was slowing down first, and then it gradually accelerated to satisfy upper-bound constraint. }
\label{T_traj_conf}
\end{figure}

We show another interesting application of our proposal. If drone and obstacle are entering in a constrained corridor, apart from putting lower bound $\mathbf{c_{min}}$, we can also put upper bound($\mathbf{c_{max}}$) in such tight spaces. Putting upper bound will ensure that drone is not slowing down too much. Sub-constraint 1 of equation \ref{new_obs_avoidance_constraint} will take a form like $\mathbf{\Upsilon_{min} \leq \Upsilon^j_{t_i} \leq \Upsilon_{max}}$. Through demonstration, we advocate the usage of upper bound constraint in tight spaces. Upper bound ensures that drone is within certain range of obstacle, which will reduce deviation in drone trajectory, thus ensuring no collision with surrounding walls. We consider a case where drone and obstacle are entering in a constrained corridor at same time. In this case, we take $\mathbf{c_{min}=30\%}$ and $\mathbf{c_{max}=60\%}$. In other words at least $\mathbf{30\%}$ confidence contours of drone and obstacle can not penetrate into each other, while $\mathbf{60\%}$ confidence contours can't have $\mathbf{0}$ overlap at any time instant during the trajectory. Our planning horizon is of $\mathbf{40}$ timesteps, each of duration $\mathbf{\tau = 0.3}$ seconds. In this case, we take  $ \mathbf{\Sigma_{drone}} =
\left(\protect\begin{smallmatrix}
      0.02 & 0.01 & 0 \\
      0.01 & 0.02 & 0 \\
      0 & 0 & 0.02 \\
      \protect\end{smallmatrix}\right)$ and $
\mathbf{\Sigma_{obs}} =
\left(\protect\begin{smallmatrix}
      0.03 & 0.02 & 0 \\
      0.02 & 0.03 & 0 \\
      0 & 0 & 0.02 \\
      \protect\end{smallmatrix}\right)$. Before solving an MPC, these matrices are scaled up to incorporate radius of drone and obstacle.  We consider non-isotropic uncertainty for this demonstration and show efficacy of our algorithm under tight spaces.  In figure \ref{T_traj_results}, We show snippets of various time-instances. The walls are modeled as stationary obstacles. For example, in figure \ref{T_traj_results}(a), drone and obstacle are entering in the corridor. Both upper bound and lower bound constraints are enforced and we can see that drone is able to maintain sufficient distance from the wall as well as the obstacle while respecting the constraints. An absence of upper-bound constraint results in slowing down of the drone and we encounter a longer time for flight completion. Having upper bound favorably changes the velocity profile to complete trajectory in faster time. This shows usefulness of our proposal in tightly bounded spaces. We can use such modeling in object tracking/following. In figure \ref{T_traj_conf}, we show confidence plot for entire trajectory.

\section{Conclusion and Future work}
\label{Conclusion and Future work_section}
In this paper, a novel approach to dynamic collision avoidance under uncertainty in state of robot and obstacles, modeled using Gaussian distribution, has been proposed. It has been derived by using area of overlap of the Gaussian distributions, which has unique characterization for a given confidence contour. The proposed algorithm, integrated with Linear MPC, optimizes over velocity profile and overlap parameter space to generate a navigation path in constrained dynamic environment. The paper puts forward results for two safety critical configurations: Antipodal configuration and Constrained Corridor setting. The findings of this study have been validated for various other possible scenarios using numerical simulations. In future, we intend to model actuation dynamics into overlap of Gaussian framework and attempt to solve for challenging scenarios with unbounded covariances. 




\bibliographystyle{IEEEtran}
\bibliography{root.bbl}

\begin{thebibliography}{10}
\providecommand{\url}[1]{#1}
\csname url@rmstyle\endcsname
\providecommand{\newblock}{\relax}
\providecommand{\bibinfo}[2]{#2}
\providecommand\BIBentrySTDinterwordspacing{\spaceskip=0pt\relax}
\providecommand\BIBentryALTinterwordstretchfactor{4}
\providecommand\BIBentryALTinterwordspacing{\spaceskip=\fontdimen2\font plus
\BIBentryALTinterwordstretchfactor\fontdimen3\font minus
  \fontdimen4\font\relax}
\providecommand\BIBforeignlanguage[2]{{%
\expandafter\ifx\csname l@#1\endcsname\relax
\typeout{** WARNING: IEEEtran.bst: No hyphenation pattern has been}%
\typeout{** loaded for the language `#1'. Using the pattern for}%
\typeout{** the default language instead.}%
\else
\language=\csname l@#1\endcsname
\fi
#2}}

\bibitem{mithun_arxiv}
M.~Babu, R.~R. Theerthala, A.~K. Singh, B.~Gopalakrishnan, and K.~M. Kirshna,
  ``Model predictive control for autonomous driving considering actuator
  dynamics,'' 2018.

\bibitem{kim2011control}
D.~Kim, J.~Kang, and K.~Yi, ``Control strategy for high-speed autonomous
  driving in structured road,'' in \emph{Intelligent Transportation Systems
  (ITSC), 2011 14th International IEEE Conference on}.\hskip 1em plus 0.5em
  minus 0.4em\relax IEEE, 2011, pp. 186--191.

\bibitem{katrakazas2015real}
C.~Katrakazas, M.~Quddus, W.-H. Chen, and L.~Deka, ``Real-time motion planning
  methods for autonomous on-road driving: State-of-the-art and future research
  directions,'' \emph{Transportation Research Part C: Emerging Technologies},
  vol.~60, pp. 416--442, 2015.

\bibitem{schwarting2017parallel}
W.~Schwarting, J.~Alonso-Mora, L.~Pauli, S.~Karaman, and D.~Rus, ``Parallel
  autonomy in automated vehicles: Safe motion generation with minimal
  intervention,'' in \emph{Robotics and Automation (ICRA), 2017 IEEE
  International Conference on}.\hskip 1em plus 0.5em minus 0.4em\relax IEEE,
  2017, pp. 1928--1935.

\bibitem{garimella2017robust}
G.~Garimella, M.~Sheckells, and M.~Kobilarov, ``Robust obstacle avoidance for
  aerial platforms using adaptive model predictive control,'' in \emph{Robotics
  and Automation (ICRA), 2017 IEEE International Conference on}.\hskip 1em plus
  0.5em minus 0.4em\relax IEEE, 2017, pp. 5876--5882.

\bibitem{wiki:Mahalanobis_distance}
Wikipedia, ``{Mahalanobis distance} --- {W}ikipedia{,} the free encyclopedia,''
  \url{https://en.wikipedia.org/wiki/Mahalanobis\_distance}.

\bibitem{van2011lqg}
J.~Van Den~Berg, P.~Abbeel, and K.~Goldberg, ``Lqg-mp: Optimized path planning
  for robots with motion uncertainty and imperfect state information,''
  \emph{The International Journal of Robotics Research}, vol.~30, no.~7, pp.
  895--913, 2011.

\bibitem{gopalakrishnan2017prvo}
B.~Gopalakrishnan, A.~K. Singh, M.~Kaushik, K.~M. Krishna, and D.~Manocha,
  ``Prvo: Probabilistic reciprocal velocity obstacle for multi robot navigation
  under uncertainty,'' 2017.

\bibitem{gopalakrishnan2015closed}
B.~Gopalakrishnan, A.~K. Singh, and K.~M. Krishna, ``Closed form
  characterization of collision free velocities and confidence bounds for
  non-holonomic robots in uncertain dynamic environments,'' in
  \emph{Intelligent Robots and Systems (IROS), 2015 IEEE/RSJ International
  Conference on}.\hskip 1em plus 0.5em minus 0.4em\relax IEEE, 2015, pp.
  4961--4968.

\bibitem{singh2013reactive}
A.~K. Singh and K.~M. Krishna, ``Reactive collision avoidance for multiple
  robots by non linear time scaling,'' in \emph{Decision and Control (CDC),
  2013 IEEE 52nd Annual Conference on}.\hskip 1em plus 0.5em minus 0.4em\relax
  IEEE, 2013, pp. 952--958.

\bibitem{scp}
S.Boyd, ``sequential convex programming,''
  \url{https://web.stanford.edu/class/ee364b/lectures/seq_slides.pdf}, 2008.

\bibitem{nowakowska2014tractable}
E.~Nowakowska, J.~Koronacki, and S.~Lipovetsky, ``Tractable measure of
  component overlap for gaussian mixture models,'' \emph{arXiv preprint
  arXiv:1407.7172}, 2014.

\bibitem{anderson1962classification}
T.~W. Anderson, R.~R. Bahadur, \emph{et~al.}, ``Classification into two
  multivariate normal distributions with different covariance matrices,''
  \emph{The annals of mathematical statistics}, vol.~33, no.~2, pp. 420--431,
  1962.

\bibitem{mathematica}
W.~R. Inc., ``Mathematica, {V}ersion 11.0,'' champaign, IL, 2016.

\bibitem{cvx}
M.~Grant and S.~Boyd, ``{CVX}: Matlab software for disciplined convex
  programming, version 2.1,'' \url{http://cvxr.com/cvx}, Mar. 2014.

\bibitem{cvxopt}
M.~Andersen and L.~Vandenberghe, ``{CVXOPT}: A python package for convex
  optimization, version 1.1.9,'' \url{http://cvxopt.org/index.html}, 2016.

\bibitem{furrer2016rotors}
F.~Furrer, M.~Burri, M.~Achtelik, and R.~Siegwart, ``Rotors—a modular gazebo
  mav simulator framework,'' in \emph{Robot Operating System (ROS)}.\hskip 1em
  plus 0.5em minus 0.4em\relax Springer, 2016, pp. 595--625.

\bibitem{koenig2004design}
N.~Koenig and A.~Howard, ``Design and use paradigms for gazebo, an open-source
  multi-robot simulator,'' in \emph{Intelligent Robots and Systems, 2004.(IROS
  2004). Proceedings. 2004 IEEE/RSJ International Conference on}, vol.~3.\hskip
  1em plus 0.5em minus 0.4em\relax IEEE, 2004, pp. 2149--2154.

\bibitem{raemaekers2007design}
A.~Raemaekers, ``Design of a model predictive controller to control uavs.''

\end{thebibliography}

\end{document}